\documentclass[12pt]{article}
\usepackage{amsmath}
\usepackage{amssymb}
\begin{document}
\newtheorem{theorem}{Theorem}[section]
\newtheorem{remark}[theorem]{Remark}
\newtheorem{mtheorem}[theorem]{Main Theorem}
\newtheorem{observation}[theorem]{Observation}
\newtheorem{proposition}[theorem]{Proposition}
\newtheorem{lemma}[theorem]{Lemma}
\newtheorem{note}[theorem]{}
\newtheorem{extlemma}[theorem]{Ext-Lemma}
\newtheorem{corollary}[theorem]{Corollary}
\newtheorem{example}[theorem]{Example}
\newtheorem{definition}[theorem]{Definition}
\newtheorem{question}[theorem]{Question}
\newtheorem{btheorem}[theorem]{Black Box Theorem}
\newtheorem{step-lemma}[theorem]{Step Lemma}

\renewcommand{\labelenumi}{(\roman{enumi})}
\newcommand{\dach}[1]{\widehat{\vphantom{#1}}}
\numberwithin{equation}{section}
\def\Z{{ \mathbb Z}}
\def\E{{ \mathbb E}}
\def\N{{ \mathbb N}}
\def\BB{ \mathbb B}
\def\DD{ \mathbb D}
\def\GG{ \mathbb G}
\def\BBB{ B_\BB}
\def\R{{\bf R}}
\def\C{{\mathcal{C}}}
\def\D{{\hat{D}}}
\def\Q{{\mathbb Q}}
\def\G{\hat{G}}
\def\T{{\cal T}}
\def\h{{\mathcal H}}
\def\V{{\mathfrak V}}
\def\X{{\mathfrak X}}
\def\Y{{\mathfrak Y}}
\def\RX{R\langle X \rangle}
\def\RXa{R\langle X_\alpha \rangle}
\def\RrXa{R(X_\alpha)}
\def\RrYa{R(Y_\alpha)}

\def\RaX{R_\alpha\langle X_\alpha \rangle}
\def\RraX{R_\alpha (X_\alpha)}
\def\K{{\mathfrak K}}
\def\qa{(R_\alpha, G_\alpha, \sigma_\alpha, \sigma_{\alpha *})}
\def\qb{(R_\beta, G_\beta, \sigma_\beta, \sigma_{\beta *})}
\def\R{{\bf R}}
\def\D{\widehat{D}}
\def\A{\widehat{A}}
\def\G{\hat{G}}
\def\T{{\cal T}}
\def\B{\widehat{B}}
\def\BC{\widehat{B_C}}
\def\BCC{\widehat{B_\C}}
\def\restr{\restriction}
\def\Aut{{\rm Aut\,}}
\def\Im{{\rm Im\,}}
\def\ker{{\rm ker\,}}
\def\inf{{\rm inf\,}}
\def\sup{{\rm sup\,}}
\def\Br{{\rm Br\,}}
\def\Yphi{Y_{[\phi]}}
\def\Ypsi{Y_{[\psi]}}
\def\Xphi{X_{\tilde{\phi}}}
\def\Xpsi{X_{\tilde{\psi}}}
\def\a{\alpha}
\def\abar{\overline{\alpha}}
\def\aa{{\bf a}}
\def\to{\rightarrow}
\def\arr{\longrightarrow}
\def\sigmaa{{\bf \Sigma_a}}
\def\End{{\rm End\,}}
\def\Ines{{\rm Ines\,}}
\def\Hom{{\rm Hom\,}}
\def\restr{\upharpoonright}
\def\Ext{{\rm Ext}\,}
\def\dom{{\rm dom}\,}
\def\Hom{{\rm Hom}\,}
\def\End{{\rm End}\,}
\def\Aut{{\rm Aut}\,}
\def\ker{{\rm ker}\,}
\def\Ann{{\rm Ann}\,}
\def\defe{{\rm def}\,}
\def\rk{{\rm rk}\,}
\def\crk{{\rm crk}\,}
\def\nuc{{\rm nuc}\,}
\def\Dom{{\rm Dom}\,}
\def\Im{{\rm Im}\,}
\def\Yphi{Y_{[\phi]}}
\def\Ypsi{Y_{[\psi]}}
\def\Xphi{X_{\tilde{\phi}}}
\def\Xpsi{X_{\tilde{\psi}}}
\def\a{\alpha}
\def\abar{\overline{\alpha}}
\def\aa{{\bf a}}
\def\ra{\rightarrow}
\def\arr{\longrightarrow}
\def\iff{\Longleftrightarrow}
\def\sigmaa{{\bf \Sigma_a}}
\def\mm{{\mathfrak m}}
\def\F{{\mathfrak F}}
\def\X{{\mathfrak X}}
\def\Diam{\diamondsuit}
\def\mapdown#1{\Big\downarrow\rlap{$\vcenter{\hbox{$\scriptstyle#1$}}$}}
\def\cRp{\rm{cRep}$_2R$ }
\def\cR{c$R_2$}

\title{{\sc Generalized $E$-Rings}
\footnotetext{This work is supported by the project
No. G-0294-081.06/93 of the German-Israeli
Foundation for Scientific Research \& Development\\
AMS subject classification (2000):\\
primary: 20K20, 20K30;  \\
secondary: 16S60, 16W20;\\
Key words and phrases: indecomposable modules, $E$-rings, homotopy theory \\
GShS 681 in Shelah's list of publications} }

\author{ R\"udiger G\"obel, Saharon Shelah and Lutz Str\"ungmann}
\maketitle
\begin{abstract}
A ring $R$ is called an $E$-ring if the canonical homomorphism
from $R$ to the endomorphism ring $\End(R_\Z)$ of the additive
group $R_\Z$, taking any $r \in R$ to the endomorphism left
multiplication by $r$ turns out to be an isomorphism of rings. In
this case $R_\Z$ is called an $E$-group. Obvious examples of
$E$-rings are subrings of $\Q$. However there is a proper class
of examples constructed recently, see \cite{DMV}. $E$-rings come
up naturally in various topics of algebra, see the introduction.
So its not surprising that they were investigated thoroughly in
the last decade, see \cite{DG, Pi, CRT, Fa, Li}. This also led to
a generalization: an abelian group $G$ is an $\E$-group if there
is an epimorphism from $G$ onto the additive group of  $\End(G)$.
If $G$ is torsion-free of finite rank, then $G$ is an $E$-group if
and only if it is an $\E$-group, see \cite{FHR}. The obvious
question was raised a few years ago which we will answer by
showing that the two notions do not coincide. We will apply
combinatorial machinery to non-commutative rings to produce an
abelian group $G$
 with (non-commutative) $\End(G)$ and the desired epimorphism
with prescribed kernel $H$. Hence, if we let $H=0$, we obtain a
non-commutative ring $R$ such that $\End(R_{\Z}) \cong R$ but $R$
is not an $E$-ring.

\end{abstract}

\section{Introduction}
Easy examples for $E$-rings are subrings of $\Q$ and further
examples up to size $2^{\aleph_0}$ coming from $p$-adic integers
are known for a long time. Details on those rings can be found in
the two monographs \cite{Fe1, Fe2}. Arbitrarily large $E$-rings
are constructed more recently first in \cite{DMV} and then in
\cite{DG}. The examples in \cite{DG} share additional properties
related with their automorphism groups. $E$-rings are important
in connection with (strongly) indecomposable modules, see
\cite{APRVW}, Pierce \cite{Pi} and Paras \cite{Pa}. They are also
crucial for investigating problems in homotopy theory, see the
work of Farjoun \cite{Fa}, \cite{Li} and Casacuberta, Rodriguez,
Tai \cite{CRT} or \cite{Ro}. So it is not surprizing that they
were studied in detail over the last decade. In this connection,
some years ago Feigelstock, Hausen and Raphael extended the
notion of $E$-rings and called an abelian group $G$ an
($EE$-group, we use the shortcut) $\E$-group if we find an
epimorphism from $G$  onto the additive group of $\End(G)$. The
results are published only recently in \cite{FHR}. In particular
the authors show that, if $G$ is torsion-free of finite rank,
then $G$ is an $E$-group if and only if $G$ is an $\E$-group -
the two notions coincide for groups of finite rank. The obvious
question was clear: It was asked whether this is true in general.
Here we want to answer this open question. The nice feature about
this problem is the fact that combinatorial arguments are needed
and must be applied to non-commutative ring theory to produce the
required examples.

We recall the two definitions.

\begin{definition} If $R$ is a ring, then $\delta: R \arr End(R_\Z)$ denotes the homomorphism which takes any $r \in R$ to the $\Z$-endomorphism $\delta(r)$ which is multiplication by $r$ on the left. If this homomorphism is an isomorphism, then $R$ is called an $E$-ring and $R_\Z$ is called an $E$-group.
\end{definition}
Recall from Schultz \cite{Sch} that $E$-rings are necessarily commutative. The rings $R$ we will construct are definitely not commutative and so not $E$-rings. However they are close to $E$-rings in the sense that the following definition holds.

\begin{definition} If $R$ is the endomorphism ring of
some abelian group $G$ and there is an epimorphism $G \arr R \arr
0$ with kernel $H$, then $G$ is called an $\E(H)$-group.
Moreover, $G$ is called an $\E$-group if $G$ is an $\E(H)$-group
for some abelian group $H$ and $G$ is a strong $\E$-group if $G$
is an $\E(0)$-group, i.e. the epimorphism is an ismorphism.
\end{definition}

Note that trivially a strong $\E$-group $G$ is in fact a ring $R$
with $R_{\Z}=G$ and such that $\End(R_{\Z}) \cong R$ as rings.

One of the main results in a paper by Feigelstock, Hausen and Raphael \cite{FHR} is the following
\begin{theorem} Let $G$ be any torsion-free abelian group of finite rank. Then $G$ is an $E$-group if and only if $G$ is an $\E$-group.
\end{theorem}

We want to prove the following result which complements the theorem by Feigelstock, Hausen and Raphael and answers their problem.

\begin{theorem}
For any infinite cardinal $\lambda = \mu^+$ with
$\mu^{\aleph_0}=\mu$ we will find an $\aleph_1$-free abelian
group $G$ of cardinality $|G| = \lambda$ which is a (strong)
$\E$-group.
\end{theorem}

\begin{remark}
A modification of our construction would also ensure that the
constructed $\E(H)$-groups are proper in the sense that they are
not strong $\E$-groups. All one has to do is to satisfy that
$G_{\alpha}/H \not\cong G_{\alpha}$ for all $\alpha < \lambda$
where $G=\bigcup\limits_{\alpha < \lambda}G_{\alpha}$ is the
$\E(H)$-group. Then the group $G$ can not be a strong $\E$-group.
\end{remark}

An $R$-module is $\aleph_1$-free if all its countable subgroups
are free. Recall that endomorphism rings of $\aleph_1$-free
abelian groups have $\aleph_1$-free additive group. The key tool
of this paper can be found in Section 2, where we invesitgate
non-commutative polynomial rings over rings. The construction of
$G$ is based on a strong version of (Shelah's) Black Box as stated
in \cite{Sh} or slightly modified in \cite{GW}, see also
\cite{CG}. The paper is based on our notes from 1998 which just
haven't been in final form for publication.

\section{Almost free rings over non-commuting variables}

Polynomial rings $R[X]$ over a ring $R$ with commuting variable $X$
are obviously free $R$-modules. We will extend this to the non-commutative
case, showing that the non-commutative polynomial ring
$\RX$ as abelian group is $\aleph_1$-free if $R_\Z$ is so.
This will be needed in Section 3. The ring construction is easy and
well-known, see Bourbaki \cite[pp. 216, 446 ff.]{Bo}. Let $R$ be a ring of
characteristic $0$, then $M$ denotes all monomials, the elements
\begin{eqnarray}\label{monom}
           r_1 X^{i_1} \cdots r_n X^{i_n}, \
 i_j \neq 0 \ (j < n), \  r_j \in R \setminus \{1\} \  (1 < j \leq n)
\mbox{ and } r_j \neq 0, (j \leq n).
\end{eqnarray}

The ring $\RX$ is formally generated by all sums of the monomials
in $M$. First we assume that $R_{\Z} = \bigoplus\limits_{t \in T}
\Z t$ is freely generated as abelian group by $T = \{t_i \;\;: i
\in I\}$ and $t_0 = 1$ without loss of generality. Then the
multiplication on $R_{\Z}$ is coded into the so-called
``constants of structure'', see Bourbaki \cite [p. 437] {Bo},
which are

\begin{eqnarray} \label{struc}
          \gamma^i_{jk} \in \Z \ (i,j,k \in \Z)
        \mbox{ with } t_j \cdot t_k = \sum_{i \in I} \gamma^i_{jk}t_i
\end{eqnarray}

The constants are well-defined by independence. We want to use $T$ to find a
representation for elements in $\RX$. Any element in $R$ can be expressed as
a sum over $T$ with coefficients in $\Z$. If $r \in \RX$ is a sum of monomials
in $M$, then we may substitute any $r_i \in R$ from (\ref{monom}) as
indicated, and the polynominal in $M$ turns into a sum of what we will call
$T$-monomials.

\begin{eqnarray}\label{Tmonom}
m = t_{i_1} X^{j_1} \cdots t_{i_n} X^{j_n}\in TM, \mbox{ with }
 j_k \neq 0 \ (k < n), \  t_{i_k} \neq 1  \  (1 < k \leq n).
\end{eqnarray}

Hence $\RX$ is generated as abelian group by all $T$-monomials and
multiplication is ruled by the structure constants (\ref{struc}) restricted
to $T$-monomials. Hence it seems plausible that the next lemma holds. We use
the structure constants on $R$ to define the ring multiplication on $\RX$. If
$m, m' \in TM$, let $m = \overline{m} t, m' = t' \overline{m'}$ with $\overline{m}, \overline{m'}
\in TM$ and $\overline{m}$ ending with $X^{\tilde {n}}, \overline{m'}$ beginning with
$X^{\tilde{m'}}$. Then the product is defined by cases.

\begin{eqnarray}\label{product}
 \mbox{ If } t \neq 1 = t', \mbox{ then } mm' = \overline m t m',
 \mbox{ and if } t' \neq 1 = t, \mbox{ then } mm' = \overline m t' m', \\
 \mbox{ if } t' = 1 = t, \mbox{ then } mm' = \overline m m',
 \mbox{ with `middle term' } X^{\tilde m + \tilde {m'}}
\end{eqnarray}
\begin{eqnarray}\label{product1}
 \mbox{ if } t' \neq 1 \neq t, \mbox{ (\ref{struc}) for }
 tt' = \sum_{i\in I} \gamma^i_{jk} t_i, t = t_j, t' = t_k,
 \mbox{ is } mm' = \sum_{i\in I} \gamma^i_{jk} m t_i \overline{m'}.
\end{eqnarray}

If $t_0 = 1$ is involved in the last equation, then the summand $\gamma^0_{jk}
m \overline{m'}$ has a `middle term' $X^{\tilde{m} + \tilde{m}}$ as in case
$t = t'= 1$. So clearly $m m'$ is a sum of $T$-monomials.

\begin{lemma}\label{free}
If $R$ is a ring with $R_\Z = \bigoplus\limits_{t \in T} \Z t,  \
T = \{ t_i : i \in I \}$ and structure constants (\ref{struc}) as
above, then let $R'_\Z = \bigoplus\limits_{m \in TM}\Z m $ be the
direct sum taken over all $T$-monomials $TM$ as above. The
multiplication on $R'$ is now defined by (\ref{product}) -
(\ref{product1}) and it follows that $R' = \RX$ as rings with
center ${\mathfrak z}R' = 1\Z$.
\end{lemma}
{\bf Proof:}
Reducing elements in $\RX$ to sums of $T$-monomials we have seen
that $R' = \RX$ as sets. Moreover, considering sums of $T$-monomials it is
obvious that $R'_{\Z} = \RX_{\Z}$ as abelian groups. We finally must show
equality as rings. There are two natural ways to do this. Either we extend
the structure constants on $R$ to $R', \RX$ and check their ring properties
(see Bourbaki \cite [p.438] {Bo}) or we check that $R'$-multiplication is
$\RX$-multiplication. We prefer the second way. Hence we must show that
(\ref{struc}), (\ref{product}) - (\ref{product1}) and its linear extension
define uniquely the ring structure on $R'$ which is the same as $\RX$.
Any case of the ring laws reduces to consider distributivity of a product
\begin{eqnarray}\label{distr}
(zm + z''m'') m'
\end{eqnarray}
with $z,z'' \in \Z, \ m, m' $ as in (\ref{product}) - (\ref{product1}) and
$m'' = \overline{m''}t''$ similar to $m$ above.
If $t = t'' = 1$ or $t' = 1$,
then (\ref{distr}) becomes immediately the unique $T$-monomial $z m m' +
z'' m'' m'$ by (\ref{product}) - (\ref{product1}) and linear extension.
If $\overline{m} = \overline{m''}$, then (\ref{distr}) can be treated with (\ref{product})
and distributivity in $T$ given by (\ref{struc}):
$$\overline{m} (z t + z'' t'') t' \overline{m'} = \overline{m} (z t t' + z'' t'' t') \overline{m'}
= z \overline{m} t t' \overline{m'} + z'' \overline{m''} t''t'
\overline{m'} = z m m' + z'' m'' m'$$ \noindent showing the
unique $\RX$-multiplication (by $m'$) in this case. If
$\overline{m} \neq \overline{m'}$, then $z m m'$ and $z'' m''m'$
are sums of {\bf distinct} $T$-monomials, hence (\ref{distr})
defines the unique $\RX$-multiplication by linear extension.
$\hfill{\square}$
\medskip

Our main interest in Lemma \ref{free} is extracted as the following
\begin{corollary}\label{freeR}
(a)  If $R$ is the ring above with $R_\Z$ free, then
$\RX_\Z /R_\Z$ is free. \\
(b) If $R_\Z$ is $\aleph_1$-free, then $\RX_\Z /R_\Z$ is
$\aleph_1$-free.
\end{corollary}
{\bf Proof:} By Lemma \ref{free} we have $T\leq TM$ and $TM$ is a $\Z$-basis
of $\RX_\Z$, hence (a) follows. For (b) choose any countable set $C$ of $R$
and let $R_C = \langle\langle C \rangle\rangle$ be the subring generated by
$C$ and $R_C \langle X \rangle$ the subring of $\RX$ generated by $C$ and
$X$ respectively. The ring $R_C$ is free by hypothesis and $R_C \langle X
\rangle / R_C$ is free by (a), hence (b) follows.
$\hfill{\square}$

If $J$ denotes the sum of all those monomials with at least one factor $X$,
then $J$ is a two sided ideal of $\RX$ (of 'non-constant' polynomials). We
have $\RX = R \oplus J$ and if $R$ is a field, then $J$ is a maximal ideal of
$\RX$. Separating summands of higher order becomes more complicated and
fortunately is not needed.

\begin{corollary}\label{Rsum}
Let $R$ be as above. The set $J$ of sums of monomials which are not constant
are a two-sided ideal of $\RX$ and $\RX = R \oplus J$ is a ring split-extension.
\end{corollary}

\section{A class of quadruples for constructing $\E$-groups}
We want to find an abelian group $G$ with $R = \End G$ and
$\sigma: G \arr R$ an epimorphism with prescribed kernel
$\ker(\sigma)=H$ and pure image $\Im(\sigma) \subseteq_* R$. Hence
$G$ is a left $R$-module and the epimorphism $\sigma$ induces a
$\Z$-homomorphism
$$ \sigma_*: G \arr R $$
such that $\sigma_*(x) \in R = \End G$ is defined by
\begin{eqnarray}\label{sigma}
\sigma_*(x)(y) = \sigma(y)x \mbox{ for all } y \in G.
\end{eqnarray}
The construction of $(G,R,\sigma, \sigma_*)$ is inductively
extending approximations of such quadruples such that the final
one is as required. Let $H$ be a fixed but arbitrary abelian
group which is $\aleph_1$-free. We say that
\begin{definition} \label{quad}
The quadruple $q = (R_q,G_q,\sigma_q, \sigma_{q*}) = (R,G,\sigma,\sigma_*)$
belongs to the class $\K_H$ if the following holds.\\
(a) \ \ $R$ is a unital ring of characteristic $0$.\\
(b) \ \ $G$ is a left $R$-module, \\
(c) \ \ $G_\Z, R_\Z$ are $\aleph_1$-free, \\
(d) \ \ $\sigma: G \arr R$ is a $\Z$-homomorphism with $\ker \sigma = H$, and $1_R \in \Im(\sigma) \subseteq_* R$\\
(e) \ \ $\sigma_*: G \arr \End G$ is a $\Z$-homomorphism defined by
      (\ref{sigma}) and $\sigma_*(g) \in R$ for each $g \in G$,
      hence $\sigma_*: G \arr R$,\\
(f) \ \ $\Ann_RG = 0$ and $R \subseteq_* \End(G)$.
\end{definition}
We often use $q_\alpha = \qa \in \K_H$ for those quadruples. The
next lemma is used to prove Proposition \ref{K}.
\begin{lemma}\label{endofree}
(a)  If $G$ is a free abelian group of finite rank, then $\End G$ is a
free abelian group.\\
(b)  If $G$ is $\aleph_1$-free, then $\End G$ is an $\aleph_1$-free abelian
group as well.
\end{lemma}
{\bf Proof:} Any endomorphism of $G$ can be represented as an element of
the cartesian product $G^G$
and $\End G \subseteq G^G$ as abelian group. However $\aleph_1$-freeness is
closed under cartesian products and subgroups, hence $\End G$ is $\aleph_1$-free.
If $G$ is freely generated by a finite set $E$, we may replace $G^G$ by
$G^E$ which is free and (a) follows.
$\hfill{\square}$

\medskip

The class $\K_H$ of quadruples is partially ordered by inclusion,
i.e. if
$$q = (R, G, \sigma, \sigma_*), q' = (R', G', \sigma', \sigma_*') \in \K_H,$$
then
$$q \leq q' \mbox{ if and only if } R \subseteq R', G \subseteq G', \sigma
\subseteq \sigma'\mbox{ and (hence) } \sigma_* \subseteq
\sigma_*'.$$ The following proposition will ensure that our
construction of an $\E(H)$-group will take place within $\K_H$.

\begin{proposition}\label{K}
(a)  If $q_i = (R_i, G_i, \sigma_i, \sigma_{i*}), \ (i \in I)$ is
a continuous, ascending chain of quadruples in $\K_H$, and $G_{i}
\subseteq_* G_{i+1}$, $1_{R_{i}}=1_{R_{i+1}}$ for all $i \in I$,
then
$$ \bigcup\limits_{i \in I} q_i = (\bigcup\limits_{i \in I}R_i,
\bigcup\limits_{i \in I}G_i, \bigcup\limits_{i \in I}\sigma_i,
\bigcup\limits_{i \in I} \sigma_{i*}) \in \K_H.$$
(b)  $(\Z, \Z \oplus H, \sigma, \sigma_*) \in \K_{H}$ where $\sigma: \Z \oplus H \rightarrow \Z$ is the canonical projection.\\
(c)  If $q \in \K_H$, then there exists $q \leq q' \in \K_H$ such
that $R_q \subseteq \sigma_{q'}(G_{q'})$ and $G_q \subseteq_*
G_{q'}$.
\end{proposition}
{\bf Proof:} (a) By continuity and (\ref{sigma}) we have
$\sigma_*=(\bigcup\limits_{i \in I} \sigma_i)_* =
\bigcup\limits_{i \in I} \sigma_{i*},$ hence $\sigma_*: G
\rightarrow R$ with pure image, where $G=\bigcup\limits_{i \in
I}G_i$ and $R=\bigcup\limits_{i \in I}R_i$. Moreover, $R$ is a
unital ring with $1_R=1_{R_i}$ for all $i \in I$ and $G$ is a left
$R$-module. $G_{\Z}$ is $\aleph_1$-free since $G_i$ is pure in
$G_{i+1}$ for all $i \in I$ and therefore $R \subseteq
\End_{\Z}G$ is $\aleph_1$-free as well by Lemma \ref{endofree}
(b). Note that $\Ann_RG=0$. Finally, $\ker\sigma=H$ and $1_R \in
\Im(\sigma)$ is clear. Hence it remains to show that $R$ is pure
in $\End_{\Z}G$. Assume that $\varphi \in \End_{\Z}G$ and
$n\varphi \in R$ for some integer $n$. Thus there is some $i \in
I$ such that $n \varphi=r' \in R_i \subseteq_* \End_{\Z}G_i$. We
claim that $\varphi \restriction_{G_j} \in \End_{\Z}G_j$ for all
$j \geq i$. Let $g \in G_j$, then $n \varphi(g)=r'g \in G_j$ and
hence $\varphi(g) \in G_j$ by purity. Thus $\varphi
\restriction_{G_j}=r_j \in R_j$ for all $i \leq j \in I$ because
$R_j \subseteq_* \End_{\Z}G_j$. Since $n\varphi=r'$, we obtain
$nr_j=r'$ for all $j$ and therefore $r:=r_j=r_k$ for all
$j,k \geq i$ by torsion-freeness. Thus $\varphi=r \in R$.  \\
(b) is obvious, and\\
(c) needs some work: Let $R'' = Rx_0 \oplus Rx$ be a ring direct
sum, $x, x_0$ two central orthogonal idempotents. We put $Rx_0 =
R$ if there is no ambiguity. Let $G' = G \oplus Re$ an
$R''$-module with $R''$ acting component-wise, hence $\Ann_{R''}e
= Rx_0 = \Ann x$ and $\Ann_{R''} G = R x$. Next we extend
$\sigma$ and let $\sigma \subseteq \sigma'$ such that $\sigma'
(re) = r x$ for all $r \in R$ (and $\sigma'(g) = \sigma (g) x_0,
g \in G)$. Hence $\sigma' : G' \arr R''$ satisfies $\ker \sigma'
=\ker \sigma =H$ and $R \subseteq \Im\sigma'$ as required.
Clearly, $\Im\sigma'=\Im\sigma x_0 \oplus Rx \subseteq_* R''$
since $\Im \sigma$ was pure in $R$. Note that $\sigma'_*$ is
actually defined by (\ref{sigma}), such that $\sigma'_* : G' \arr
\End_{\Z} G'$, hence $R''$ must be enlarged for $\sigma'_* : G'
\arr R'$ in Definition \ref{quad}(e). The ring $R''$ acts by
scalar multiplication (on the left) on $G'$. The action is
faithful by hypothesis (Definition \ref{quad} (f)), hence $R''$
can be viewed as a subring of $\End_{\Z} G'$. Let $R' = R'' [\Im
\sigma'_*]_*$ be the pure unital subring of $\End_{\Z} G'$
generated by $R''$ and $\Im \sigma'_*$.\\

Obviously $q \leq q'$ for $q' = (R', G', \sigma', \sigma'_*)$ and
$\sigma' = \sigma_q'$ satisfies (c) of Proposition \ref{K}. It
remains to show that $q' \in \K_H$. We must check Definition
\ref{quad} (c), (d) and (f). We have $\ker \sigma' = \ker \sigma
= H$. Moreover, if $g \in G'$, then $\sigma'_* (g) (G') = 0$
implies $\sigma'_* (g) = 0$, hence $\Ann_{R'} G' = 0$ and (f)
follows. It remains to show that $\Im\sigma'$ is pure in $R'$. As
shown above, $\Im \sigma'$ is pure in $R''$, hence it suffices to
show that $R''$ is pure in $R'$. We will even show that $R''$ is
pure in $\End_{\Z}G'$. Let $\varphi \in \End_{\Z}G'$ and assume
that $n\varphi=rx_0 \oplus r'x \in R''$ for some integer $n$. By
the torsion-freeness of $G'$ and $R_{\Z}$ it follows that
$\varphi$ is of the form $\varphi' \oplus r^*$ for some $\varphi'
\in \End_{\Z}G$ and $r^* \in R$. Hence $r'=nr^*$ and
$n\varphi'=r$ and therefore $\varphi'=r'' \in R$ since $R$ is
pure in $\End_{\Z}G$. Thus $\varphi=r''x_0 \oplus r^*x \in R''$
and the purity of $R''$ is established. Finally, the abelian
groups $G,R_{\Z}$ are $\aleph_1$-free by hypothesis on $q$, hence
$G'$ is $\aleph_1$-free and $\End_{\Z} G'$ is $\aleph_1$- free by
Lemma \ref{endofree}, and therefore $q' \in \K_H$.
$\hfill{\square}$\\

\bigskip
\noindent
\relax From the last proof we extract a useful
\begin{definition}\label{RE}
If $\sigma_*: G \arr \End_\Z G$ is a $\Z$-homomorphism and $R$ acts faithful
on the left $R$-module $G$ by scalar multiplication, then we denote by
$$R_{\sigma_*} \subseteq_* \End_\Z G$$
the pure unital subring of $\End_\Z G$ generated by $R \subseteq
\End_\Z G$ and $\Im \sigma_*$.
\end{definition}
The next proposition provides the link to our construction in
Section 4.
\begin{proposition}\label{trans}
Let $q = (R,G,\sigma, \sigma_*) \in \K_H$. Then there exists a
`transcendental extension' $q \leq q' \in \K_H$ such that the
following holds for
$q' = (R', G', \sigma', \sigma_*')$.\\
(a) \ \ $G' = G \oplus \RX e$ as $\RX$-module where $X$ acts as identity on $G$.\\
(b) \ \ $\sigma'(g)=\sigma(g)Xe$, $\sigma'(re) = re$ if $r \in R$ and $\sigma'(re)=rXe$ else.\\
(c) \ \ $R' =  (\RX)_{\sigma'_*}$ with $\RX$ from Section 2.\\
(d) \ \ $R \subseteq \Im(\sigma)$.
\end{proposition}

{\bf Proof:} We must show that $q' \in \K_H$. The ring $R$ is
$\aleph_1$-free by hypothesis, hence $\RX$ is $\aleph_1$-free by
Corollary \ref{freeR}(b) and $R'$ is $\aleph_1$-free by Lemma
\ref{endofree}. Clearly, $\ker \sigma' =\ker \sigma =H$ and $\Im
\sigma' \subseteq_* R'$ with $R \subseteq \Im(\sigma)$ follow as
in the proof of Proposition \ref{K}. Moreover, $G'$ is
$\aleph_1$-free and $\Ann_{R'}G' = 0$,
so the proposition follows. $\hfill{\square}$\\
\noindent

We have an immediate corollary-definition.\\

\begin{corollary} If
$K_H^o = \{ q \in \K_H, \sigma_q \mbox{ maps onto } R_q \}, $ then
$\K_H^o$ is dense in $\K_H$.
\end{corollary}
{\bf Proof:} Apply Proposition \ref{K}(c) $\omega $ times and
note that Proposition \ref{K}(a) can be used because the union of
the countable sequences of rings and modules respectively are
$\aleph_1$-free by construction. Hence any $q \in \K_H$ is below
some $q' \in \K_H^o$. $\hfill{\square}$

\section{Construction of $\E$-groups by black box arguments}
The combinatorial ideas of Sections 4 and 5 can be found in
Shelah \cite{Sh}, see also the appendix {\it Shelah's Black Box}
in Corner, G\"obel \cite{CG} and this Black Box could be used.
However, we will use a stronger version of the Black Box as
developed in \cite{GW}. The Black Box needs, as usually, some
harmless alterations, which are obvious and the proof is left to
the reader.

Let $\lambda$ be some infinite cardinal, and $\{X_\alpha : \alpha
< \lambda \}$ a sequence of transcendental elements which will be
used to define ring extensions:

First we define ring extensions `locally' and let $R_{\alpha+1}
\supseteq R_{\alpha} \left(X_{\alpha} \right) = R_{\alpha}\left<
X_{\alpha} \right>=R_{\alpha}c_{\alpha} \oplus
J_{\alpha}c_{\alpha}'$ a ring direct sum with central idempotents
$c_{\alpha}$ and $c_{\alpha}'$ where $J_{\alpha}$ is defined as
in Corollary \ref{Rsum}. Since there is no danger of confusion we
usually will omit the 'place holders' $c_{\alpha}$ and
$c_{\alpha}'$. The sequence of rings $R_{\alpha}$ with
$|R_{\alpha}| < \lambda$ $(\alpha < \lambda$) will be completed
during the construction of the abelian group $G$ (with
$\End_{\Z}G = \bigcup\limits_{\alpha < \lambda}R_{\alpha}$),
taking unions at limit steps,
i.e. $R_{\alpha}=\bigcup\limits_{\beta < \alpha}R_{\beta}$ if $\alpha$ is a limit ordinal.
Note that $\bigcup\limits_{\alpha < \lambda}R_{\alpha}=\bigcup\limits_{\alpha < \lambda}R_{\alpha}\left(X_{\alpha}\right)$.\\

\noindent
Similarly we define
$$ B = \bigoplus\limits_{\alpha < \lambda} R_{\alpha} \left( X_{\alpha} \right) e_\alpha \subseteq G
\subseteq \B$$ where $\widehat{B}$ denotes the $p$-adic completion
of $B$ (as an abelian group) for some fixed prime $p$. Since we
want to apply the Black Box later on we need a free basis-module
inside $B$. Therefore, assume that our rings $R_{\alpha}$ (and
hence $R_{\alpha}(X_{\alpha}))$ ($\alpha < \lambda$) are
$\aleph_1$-free, hence homogeneous of type $\Z$. By a well-known
result \cite[Theorem 128]{G} there exists for each $\alpha <
\lambda$ a completely decomposable group $F_{\alpha} \subseteq_*
R_{\alpha}(X_{\alpha})$ such that $|R_{\alpha}(X_{\alpha})| \leq
|F_{\alpha}|^{\aleph_0}$ and
\[ F_{\alpha} \subseteq_* R_{\alpha}(X_{\alpha}) \subseteq_*
\widehat{F_{\alpha}} \] where $\widehat{F_{\alpha}}$ is the
$p$-adic completion of $F_{\alpha}$ ($\alpha < \lambda$). Note
that $F_{\alpha}$ is a free abelian group since
$R_{\alpha}(X_{\alpha})$ is $\aleph_1$-free. We collect all those
$F_{\alpha}$ ($\alpha < \lambda$) and define $F:=
\bigoplus\limits_{\alpha < \lambda}F_{\alpha}e_{\alpha}$. Thus $F$
is a free abelian group of cardinality at most $\lambda$ such that
\[ F \subseteq_* B \subseteq_* \widehat{F}. \]
Let $F_{\alpha}=\bigoplus\limits_{\varepsilon < \rho} \Z
a_{\varepsilon}$ where $\rho =|F|$. Writing $e_{(\varepsilon,
\alpha)}$ for $a_{\varepsilon}e_{\alpha}$ it follows that
$B'=\bigoplus\limits_{(\varepsilon, \alpha) \in \rho \times
\lambda} \Z e_{(\varepsilon, \alpha)}$ satisfies
$\B=\widehat{B'}$. For later use we put the lexicographic
ordering on $\rho \times \lambda$; since $\rho, \lambda$ are
ordinals $\rho \times \lambda$ is well
ordered.\\
We are ready to define supports of elements in $\hat{B}$.
If $0 \neq g \in \B$ then we can write $g = \sum\limits_{\alpha
\in I} g_\alpha e_\alpha$ and $I \subseteq \lambda, |I| \leq
\aleph_0, g_\alpha \in R_\alpha (X_\alpha)$. Moreover, each
$g_\alpha = g'_\alpha + g''_\alpha$ with $g'_\alpha \in R_\alpha$
and $g''_\alpha \in J_\alpha$.
We define:\\
\begin{eqnarray} \label{support}
\mbox{ The $\lambda$-support of } g \mbox{ is the set }
[g]_{\lambda} = \{ \alpha  < \lambda : g_\alpha \neq 0 \}. \\
\end{eqnarray}

The notion of $\lambda$-support naturally extends to subsets of
$\widehat{B}$, see again \cite{GW}.\\
On the other hand, any element $0 \not= g \in \B=\widehat{B'}$ can
be written as
\[ g=(g_{(\varepsilon, \alpha)}e_{(\varepsilon,
\alpha)})_{(\varepsilon, \alpha) \in \rho \times \lambda} \in
\widehat{B'} \subseteq \prod\limits_{(\varepsilon, \alpha) \in
\rho \times \lambda}\widehat{Z}e_{(\varepsilon, \alpha)} \] and we
define the 'usual' support of $g$ by $[g]=\{ (\varepsilon, \alpha)
\in \rho \times \lambda | g_{(\varepsilon, \alpha)} \not=0 \}$.
Note that $| [g]| \leq \aleph_0$ and that the $\lambda$-support
of $g$ is $[g]_{\lambda}=\{ \alpha < \lambda | \exists
\varepsilon < \rho : (\varepsilon, \alpha) \in [g] \}$. As usual
we may define a norm on $\widehat{B'}$ by $|| \alpha||=\alpha +1$
($\alpha < \lambda$), $||M||=\sup_{\alpha \in M}||\alpha||$ ($M
\subseteq \lambda)$ and $||g||=||[g]_{\lambda}||$ $(g \in
\widehat{B'}$), i.e. $||g||=\min \{\beta < \lambda |
[g]_{\lambda} \subseteq \beta \}$. Note, $[g]_{\lambda} \subseteq
\beta$ holds if and only if $g \in \widehat{B_{\beta}'}$ where
$B_{\beta}'=\bigoplus\limits_ {\alpha < \beta}
R_{\alpha}\left(X_{\alpha} \right)e_{\alpha}$. Finally, we also
have a ring support and ring norm defined by
$[g]_{ring}=\{\varepsilon < \rho | \exists \alpha < \lambda :
(\varepsilon, \alpha) \in [g] \}$ and
$||g||_{ring}=||[g]_{ring}||$.  \\
 We will now state a suitable version of the
Strong Black Box as developed in \cite{GW}. The proof is almost
identical with the one in \cite[Section 1]{GW} and will therefore
be left to the reader but we will
give all the necessary definitions and results.\\
Fix cardinals $\kappa \geq \aleph_0$, $\mu=\mu^{\kappa}$ such that
$\lambda=\mu^+$. We need to say what we mean by a canonical
homomorphism. For this we fix bijections $g_{\gamma}: \mu
\rightarrow \gamma$ for all $\gamma$ with $\mu \leq \gamma <
\lambda$ where we put $g_{\mu}=id_{\mu}$. For technical reasons
we also put $g_{\gamma}=g_{\mu}$ for $\gamma < \mu$. Moreover,
let $g_{(\varepsilon, \alpha)}=g_{\varepsilon} \times g_{\alpha}$
for all $(\varepsilon, \alpha) \in \rho \times \lambda$.

\begin{definition}\label{defcan1}\quad 
We define $P$ to be a {\em canonical summand} of $B'$ if
$P=\bigoplus\limits_{(\varepsilon,\alpha)\,\in\, I} \Z
e_{(\varepsilon,\alpha)}$ for some $I\subseteq\rho\times\lambda$
with $|I|\leq\kappa$ such that:
\begin{itemize}
\item
if $(\varepsilon,\alpha)\,\in I$ then
$(\varepsilon,\varepsilon)\in I$;
\item
if $(\varepsilon,\alpha)\,\in I, \alpha\in\rho$ then
$(\alpha,\varepsilon)\in I$;
\item
if $(\varepsilon,\alpha)\in I$ then
$\left(I\cap(\mu\times\mu)\right)g_{\varepsilon,\alpha}= I\cap \Im
g_{\varepsilon,\alpha}$; \ and
\item
$||P|| < \lambda^o$, where $\lambda^0=\{ \alpha < \lambda | \text{ cf}(\alpha)=\aleph_0 \}$. \end{itemize} 
Accordingly, $\varphi\!:P\rightarrow \widehat{B'}$ is said to be a
{\em canonical homomorphism} if $P$ is a canonical summand of $B'$
and $\Im\varphi\subseteq \widehat{P}$; we put $[\varphi]=[P]$,
$[\varphi]_{\lambda}=[P]_{\lambda}$ and $||\varphi||=||P||$.
\end{definition}

If we denote by $\C$ the set of all canonical homomorphisms, then
$|\C|=\lambda$ holds (see \cite{GW}). Our version of the Strong
Black Box reads as follows (compare \cite[Theorem 1.1.2.]{GW}):

\begin{btheorem}
\label{black}
Let $E \subseteq \lambda^0$ be a stationary subset of $\lambda$ with $\lambda=\mu^+$, $\mu^{\kappa}=\mu$.\\
Then there exists a family $\C^*$ of canonical homomorphisms with
the following properties:
\begin{enumerate}
\item If $\varphi \in \C^*$, then $||\varphi|| \in E$.
\item If $\varphi, \varphi'$ are two different elements of $\C^*$
of the same norm $\alpha$ then $||[\varphi]_{\lambda} \cap
[\varphi']_{\lambda}|| < ~\alpha$.
\item PREDICTION: For any homomorphism $\psi: B' \rightarrow
\widehat{B'}$ and for any subset $I$ of $\lambda$ with $|I| \leq
\kappa$ the set
\[ \{ \alpha \in E | \exists \varphi \in \C^* : ||\varphi|| =
\alpha, \varphi \subseteq \psi, I \subseteq [\varphi] \} \] is
stationary.
\end{enumerate}
\end{btheorem}

For the proof of the above Theorem we have to define an
equivalence relation on $\C$:

\begin{definition}\label{type1}\quad
Canonical homomorphism $\varphi,\,\varphi'$ are said to be {\em
equivalent} or {\em of the same type} $($notation:
 $\varphi\equiv\varphi')$, if
$[\varphi]\cap(\mu\times\mu)=[\varphi']\cap(\mu\times\mu)$ and
there exists an order-isomorphism $f\!:[\varphi]\rightarrow
[\varphi']$ such that $(x\bar{f})\varphi'=(x\varphi)\bar{f}$ for
all $x\in\dom\varphi$ where
$\bar{f}\!:\widehat{\dom\varphi}\rightarrow
\widehat{\dom\varphi'}$ is the unique extension of the
$R$-homomorphism defined by
$e_{(\varepsilon,\alpha)}\bar{f}=e_{(\varepsilon,\alpha)f}\
\left((\varepsilon,\alpha)\in [\varphi]\right)$.
\end{definition}

As in \cite{GW} it is easy to see that there are at most $\mu$
different types. Next we have to recall the definition of an
admissible sequence.

\begin{definition}
Let $\varphi_0 \subset \varphi_1 \subset \subset \cdots \subset
\varphi_n \subset \cdots (n < \omega)$ be an increasing sequence
of canonical homomorphisms.\\
Then $\left( \varphi_n \right)_{n < \omega}$ is said to be
admissible if $[\varphi_0] \cap (\mu \times \mu)=[\varphi_n]\cap
(\mu \times \mu)$ for all $n < \omega$. Also, we say that $\left(
\varphi_n \right)_{n < \omega}$ is admissible for a sequence
$\left( \beta_n \right)_{n < \omega}$ of ordinals in $\lambda$ if
$\left( \varphi_n \right)_{n < \omega}$ is admissible satisfying
$||\varphi_n|| \leq \beta_n < ||\varphi_{n+1}||$ and
$[\varphi_n]=[\varphi_{n+1}] \cap
(\beta_n \times \beta_n)$ for all $n < \omega$.\\
Moreover, two admissible sequences $\left( \varphi_n \right)_{n <
\omega}$ and $\left( \varphi_n' \right)_{n < \omega}$ are said to
be equivalent or of the same type if $\varphi_n \equiv \varphi_n'$
for all $n < \omega$.
\end{definition}

Note that the union $\bigcup\limits_{n < \omega}\varphi_n$ of an
admissible sequence $\left( \varphi_n \right)_{n < \omega}$ is
also an element of $\C$. Moreover, if we let $\T$ be the set of
all possible types of admissible sequences of canonical
homomorphisms, then clearly $|\T| \leq \mu^{\kappa}=\mu$. If
$\left( \varphi_n \right)_{n < \omega}$ is admissible of type
$\tau$ , then we also use the notion of $\tau$-admissible. As in
\cite{GW}, the following proposition is the main ingredient of
the proof of the Black Box Theorem \ref{black}.

\begin{proposition}
Let $\psi: B' \rightarrow \widehat{B'}$ be a homomorphism, $I
\subseteq \rho \times \lambda$ a set of cardinality at most
$\kappa$ and $\h=\h_{\psi,I}=\{ \varphi \in \C | \varphi
\subseteq \psi, I
\subseteq [\varphi] \}$.\\
Then there exists a type $\tau \in \T$ such that
\[ \exists \varphi_0 \in \h \forall \beta_0 \geq ||\varphi_0||
\cdots \exists \varphi_n \in \h \forall \beta_n \geq
||\varphi_n|| \cdots \] with $\left( \varphi_n \right)_{n <
\omega}$ is $\tau$-admissible.
\end{proposition}

The proof is contained in \cite{GW} and also for the proof of the
Black Box Theorem \ref{black} we refer to \cite[Section 1]{GW}.
Finally we have a corollary suitable for application.

\begin{corollary}
\label{predcor} Let the assumption be the same as in the Black Box
Theorem \ref{black}. Then there exists an ordinal $\lambda^* \geq
\lambda$ with $|\lambda^*|=\lambda$ and a family $\left(
\varphi_{\beta} \right)_{\beta < \lambda^*}$ of canonical
homomorphisms such that
\begin{enumerate}
\item $\varphi_{\beta} \in \C^*$ and $|| \varphi_{\beta}|| \in E$ for all $\beta < \lambda^*$.
\item $||\varphi_{\gamma}|| \leq ||\varphi_{\beta}||$ for all
$\gamma \leq \beta < \lambda^*$.
\item $||[\varphi_{\beta}]_{\lambda} \cap [\varphi_{\gamma}]_{\lambda}||
< ||\varphi_{\beta}||$ for all $\gamma < \beta < \lambda^*$.
\item PREDICTION: For any homomorphism $\psi: B' \rightarrow
\widehat{B'}$ and for any subset $I$ of $\lambda$ with $|I| \leq
\kappa$ the set
\[ \{ \alpha \in E | \exists \beta < \lambda^* : ||\varphi_{\beta}|| =
\alpha, \varphi_{\beta} \subseteq \psi, I \subseteq
[\varphi_{\beta}] \}
\] is stationary.
\end{enumerate}
\end{corollary}

\section{The inductive steps in the construction of $q = (R, G,
\sigma, \sigma_*)$.}

We now use induction along $\alpha < \lambda^*$ given by the Black
Box to find quadruples $q_\alpha = \qa \in \K_H$ for a fixed
$\aleph_1$-free group $H$, see Definition \ref{quad}. Let $H$ be
given and choose $q_0 \in \K_H$ arbitrary, e.g. $q_0 = (\Z,\Z
\oplus H,\sigma,\sigma_*)$ which is in $\K_{H}$ by Proposition
\ref{K} (b). Necessarily we have to assume that $|H| < \lambda$.
Let $(\varphi_{\beta})_{\beta < \lambda^*}$ be a family of
canonical homomorphisms as given by Corollary \ref{predcor}. For
any $\beta < \lambda^*$ let $P_{\beta}=\dom \varphi_{\beta}$.
Suppose that the quadruples $q_{\beta}=\qb$ are constructed for
all $\beta < \alpha$ subject
to the following conditions: \\
\begin{enumerate}
\item $R_{\beta}$ is a unital ring with
$1_{R_{\beta}}=1_{R_{\gamma}}$ for all $\gamma \leq \beta$
\item $R_\beta, G_\beta$ are $\aleph_1$-free
\item $\sigma_{\beta} (G_\beta) \subseteq R_\beta$
\item $G_{\beta} \subseteq_* G_{\beta + 1}$
\item $R_{\beta} \subseteq \Im \sigma_{\beta+1}$ if $\beta \not\in
E$
\item $G_{\beta}=\bigcup\limits_{\gamma < \beta} G_{\gamma}$ if
$\beta$ is a limit ordinal
\item $B_{\beta} = \bigoplus\limits_{\gamma < \beta}
R_{\gamma} \left(X_{\gamma} \right)e_{\gamma} \subseteq_*
G_{\beta} \subseteq_* \widehat{B_{\beta}}$ where
$\widehat{B_{\beta}}$ is the $p$-adic completion of $B_{\beta}$.
\end{enumerate}

We first have to prove a Step Lemma.

\begin{step-lemma}
\label{step} Let $P=\bigoplus\limits_{(\varepsilon, \alpha) \in
I^*} \Z e_{(\varepsilon, \alpha)}$ for some $I^* \subseteq \rho
\times \lambda^*$ and let $M$ be a subgroup of $\widehat{B'}$
with $P \subseteq_* M \subseteq_* \widehat{B'}$ which is
$\aleph_1$-free and an $R$-module, where
$R=\bigcup\limits_{\alpha \in I'}R_{\alpha}\left( X_{\alpha}
\right)$ with $I'=\{ \alpha < \lambda : \exists \varepsilon <
\rho, (\varepsilon, \alpha) \in I^* \}$. Assume that
$q=(R,M,\sigma, \sigma_*) \in \K_H$. Also suppose that there is a
set $I = \{ (\varepsilon_n,\alpha_n) : n < \omega \} \subseteq
[P]=I^*$ such that $\alpha_0 < \alpha_1 < \cdots < \alpha_n <
\cdots (n < \omega)$ and
\begin{enumerate}
\item $M=\bigcup\limits_{n < \omega}G_{\alpha_n}$,
$R=\bigcup\limits_{n < \omega}R_{\alpha_n}$;
\item $(R_{\alpha_n}, G_{\alpha_n}, \sigma_{\alpha_n}, \sigma_{\alpha_n*}) \in
\K_H$ for all $n < \omega$
\item $I_{\lambda} \cap [g]_{\lambda}$ is finite for all $g \in M$ $(I_{\lambda}=[I]_{\lambda}$).
\end{enumerate}
Moreover, let $\varphi : P \rightarrow M$ be a homomorphism which
is not multiplication by an element from $R$.\\
Then there exists an element $y \in \widehat{P}$ and an element
$q \subseteq q'=(R',M', \sigma', \sigma'_*) \in \K_H$ such that $M
\subseteq_* M' \subseteq_* \widehat{B'}$, $y \in M'$ and $y\varphi
\not\in M'$.
\end{step-lemma}

{\bf Proof:} By assumption $M$ is an $R$-module and hence the
completion $\widehat{M}$ is an $\widehat{R}$-module. Thus for any
$y \in \widehat{P}$ and $r \in \widehat{R}$ it follows that $ry$
is defined inside $\widehat{M}$ and
\[ [ry] \subseteq [y]. \]
We construct a new group $M \subseteq_* M_y$ for $y \in
\widehat{P}$ as follows: \\
Put $M_y^1=\left< M, Ry \right>_* \subseteq \widehat{B'}$. Since
$\Im \sigma$ is a pure subgroup of $R$ there is a unique extension
\[ \widehat{\sigma} : \widehat{M} \rightarrow
\widehat{R}  \] of $\sigma$ such that $\ker \widehat{\sigma} =
\ker \sigma = H$. Moreover, $\Im \widehat{\sigma}$ is pure in
$\widehat{R}$. We choose $R_y^1=R[\Im \widehat{\sigma}
\restriction_{M_y^1}]_* \subseteq \widehat{R}$ and let \[
M_y^2=\left< R_y^1M, R_y^1y \right>_* \subseteq \widehat{B'}.
\] Again we may choose $R_y^2= R_y^1[\Im \widehat{\sigma}
\restriction_{M_y^2} \subseteq \widehat{P}]_* \subseteq
\widehat{R}$. Continuing this way we obtain a sequence of groups
$M_y^n$ and rings $R_y^n$ ($n \in \omega$) such that
\[ M_y^{n+1} \text{ is an } R_y^n \text{-module and }
\widehat{\sigma} (M_y^{n}) \subseteq R_y^{n}. \] Taking
$M_y=\bigcup\limits_{n \in \omega}M_y^n$ and
$R'_y=\bigcup\limits_{n \in \omega} R_y^n$ we get that $M_y$ is an
$R'_y$-module and $\sigma_y=\widehat{\sigma} \restriction_{M_y} :
M_y \rightarrow R'_y$. It is a standard support argument to see
that $M_y$ and (hence) $R'_y$ are still $\aleph_1$-free. Finally
take $R_y=\left(R_y' \right)_{\sigma_y^*} \subseteq
\End_{\Z}(M_y)$. By construction $\Im \sigma_y$ is pure in $R_y$
and hence $q_y=(R_y, M_y, \sigma_y, \sigma_y^*) \in \K_H$.\\ At
this stage we determine $y$ more specific in order to obtain that
$\varphi \not\in \End_{\Z} M_y$. Let $y=\sum\limits_{n \in
\omega}p^ne_{(\varepsilon_n, \alpha_n)}$ and $x=\varphi(y) \in
\widehat{M}$. If $x \not\in M_y$, then choose $M'=M_y$ and
$R'=R_y$ and hence
$q_y \in \K_H$ with $\varphi \not\in \End_{\Z} M_y$.\\
 If
$x \in M_y$, then there are integers $k$ and $n$ such that
\[ p^k\varphi(y)=r_ng + r_n'y \]
for some $r_n, r_n' \in R_y^n$ and $g \in M$. It follows that
\[ \left(p^k\varphi -r_n' \right)y = r_ng. \]
Since $\left(p^k \varphi - r_n'\right) \not=0$ there is $b' \in P$
such that
\[ \left(p^k\varphi -r_n'\right)b' \not=0. \]
Note that $b'$ has finite support. Moreover, by the
cotorsion-freeness of $R$ there exists $\pi \in \widehat{R}$ such
that $\pi b \not\in M$ with $b=\left(p^k\varphi - r_n' \right)
b'$. Let $y'=y + \pi b'$. We claim that $\varphi \not\in
\End_{\Z}(M_{y'})$. By way of contradiction assume that
\[ p^l\varphi(y+\pi b') = r^*_mg^* + r_m^{'*}(\pi b'
+ y) \] for some integer $l \geq k$ and elements $r_m^*, r_m^{'*}
\in R_{y'}^n$ and $g^* \in M$. Without loss of generality, we may
assume $n=m$, hence
\[ p^l\varphi(y+\pi b') = r^*_ng^* + r_n^{'*}(\pi b'
+ y). \] Let $s=p^l/p^k$. Hence
\[ p^l\varphi(\pi b') = p^l\varphi(y + \pi b') -
sp^k\varphi(y) = r_n^*g^* + r_n^{'*}(\pi b' + y) - s(r_ng +
r_n'y) =\]
\[= (r_n^*g^* - sr_ng) + r_n^{'*}\pi b' +
(r_n^{'*} - sr_n')y. \] Since $[\pi b' ]= [b'], [\varphi(\pi b')
] = [ \varphi (b')]$ and $g^*, g \in M$ an easy support argument
shows that $r_n^{*'}=sr_n'$ and hence
\[ sp^k \varphi(\pi b') = (r_n^*g^* - sr_ng) + sr_n' \pi b' \]
and thus
\[ s\pi  (p^k \varphi(b') - r_n^{'} b')= (r_n^*g^*
- sr_ng ) \in M.\] By purity we get $\pi(p^k \varphi (b') -
r_n^{'} b' ) =\pi b \in M$ - a contradiction. Finally we put $M' =
M_{y'}$, $R' =R_{y'}$ and $q'=q_y \in \K_H$. Thus \[ M'=R'M +
\sum\limits_{k < \omega}R'y^{(k)}, \] where
$y^{(k)}=\sum\limits_{n \geq k}\frac{p^n}{p^k}e_{(\varepsilon_n,
\alpha_n)}$ or $y^{(k)}=\sum\limits_{n \geq k} \frac{p^n}{p^k}
e_{(\varepsilon_n, \alpha_n)} + \pi^{(k)}b$.
$\hfill{\square}$ \\

We will now carry on the construction to $\alpha$ and distinguish
three cases.\\

 \noindent
{\bf Case 1}: Suppose $\alpha$ is a limit ordinal. Then
$G_{\alpha}=\bigcup\limits_{\beta < \alpha}G_{\beta}$ is
$\aleph_1$-free by (iii) and hence Proposition \ref{K} (a) shows
that we can take unions, i.e. $q_{\alpha}=\bigcup\limits_{\beta <
\alpha}q_{\beta}$.
\\
\noindent {\bf Case 2}: Suppose $\alpha=\beta +1$, then
$||\varphi_{\beta} || \in \lambda^0$. Assume that $\Im
\varphi_{\beta} \not\subseteq G_{\beta}$ or $\varphi_{\beta} \in
R_{\beta}$. In this case we let $G_{\beta+1} = G_\beta \oplus
R_\beta \langle X_\beta \rangle e_\beta$ as in Proposition
\ref{trans} with $R'_{\alpha}=R_{\beta} \left< X_{\beta} \right> =
R_{\beta}\left( X_{\beta} \right)$ with $X$ acting as identity on
$G_{\beta}$. We let $R_{\alpha}=\left( R_{\alpha}'
\right)_{\sigma_{\alpha}^*}$ and by Proposition \ref{trans} it
follows that $(R_{\alpha}, G_{\alpha}, \sigma_{\alpha},
\sigma_{\alpha}^*) \in \K_H$ where $\sigma_{\alpha}$ and
$\sigma_{\alpha}^*=\sigma_{\alpha*}$ are taken from Proposition
\ref{trans}. Put $y_{\beta}=0$.
\\

\noindent {\bf Case 3}: Suppose that $\alpha=\beta + 1$ and $\Im
\varphi_{\beta} \subseteq G_{\beta}$, $\varphi_{\beta} \not\in
R_{\beta}$. In this case we try to `kill' our undesired
homomorphism $\varphi_{\beta}$ which comes from the Black Box
prediction. \\
Recall that $|| \varphi_{\beta} || \in \lambda^0$, hence there are
$(\varepsilon_n, \beta_n) \in [\varphi_{\beta} ]$ $(n \in
\omega$) such that $\beta_0 < \beta_1 < \cdots < \beta_n <
\cdots$ and $\sup_{n \in \omega}\beta_n = ||\varphi_{\beta} ||$.
Without loss of generality we may assume that $\beta_n \not\in E$
for all $n \in \omega$ and hence $G_{\beta_n +1} = G_{\beta_n}
\oplus R_{\beta_n} \left< X_{\beta_n} \right>e_{\beta_n}$. We put
$I= \{(\varepsilon_n,  \beta_n) | n < \omega \}$. Then
$I_{\lambda} \cap [g]_{\lambda}$ is finite for all $g \in
G_{\beta}$. We apply the Step Lemma ~\ref{step} to $I$ as above,
$P=\dom \varphi_{\beta}$ and $M=G_{\beta}$. Therefore there
exists an extension $q_{\alpha}=q_{\beta+1}$ of $q_{\beta}$ and
an element $y_{\beta} \in G_{\alpha}$ such that $y_{\beta}
\varphi_{\beta} \not\in G_{\alpha}$ and
$||y_{\beta}||=||\varphi_{\beta}|| = || P_{\beta}
||$. \\

Finally put
$q_H=(G_H,R_H,\sigma_H,\sigma_{H*})=\bigcup\limits_{\alpha <
\lambda}q_{\alpha} \in \K_H$. Obviously, $G_H$ has cardinality
$\lambda$ and $R \subseteq \End_{\Z}G$. Moreover,
\[ G_H = B + \sum\limits_{\beta < \lambda^*}\sum\limits_{k < \omega} R_{\beta}y_{\beta}^{(k)}. \]\\
Next we describe the elements of $G_H$.

\begin{lemma}
\label{desc1} Let $G_H$ be as above and let $g \in G_H \backslash
B$. Then there are $k < \omega$ and a finite subset $N$ of
$\lambda^*$ such that $g \in B + \sum\limits_{\beta \in
N}R_{\beta}y_{\beta}^{(k)}$ and $[g]_{\lambda} \cap
[y_{\beta}]_{\lambda}$ is infinite if and only if $\beta \in N$.
In particular, if $||g||$ is a limit ordinal then
$||g||=||y_{\beta_*}||=|| \varphi_{\beta_*}||$ for $\beta_*=\max
N$.
\end{lemma}

{\bf Proof:} Let $g\in G_H=B +
\sum\limits_{\beta<\lambda^*}\sum\limits_{n<\omega}R_{\beta}y_\beta^{(n)}$.
Then there are a finite subset $N'$ of $\lambda^*$, $b\in B,\,
k\in\omega,\, a_{\beta,n}\in R_{\beta} \ (\beta\in N',\,n\leq k)$
such that\\ 
\centerline{$g=b + \sum\limits_{\beta\in N'}\sum\limits_{n \leq
k}a_{\beta,n}y_{\beta}^{(n)}$.}\\ 
Since $y_\beta^{(n)}-\frac{p^k}{p^n}y_\beta^{(k)}\in B'\subseteq
B$ this expression reduces to
 \\ 
\centerline{$g=b' + \sum\limits_{\beta\in N'}a_\beta y_{\beta}^{(k)}$}\\ 
for some $a_\beta\in R_{\beta}\ (\beta\in N'),\, b' \in B'$.
Putting $N=\{\beta\in N'\,|\, a_\beta\not=0\}\
(N\not=\emptyset\text{ for } g \notin B)$ the conclusion of the
lemma follows since $[y_\beta]_{\lambda}\cap
[y_{\beta'}]_{\lambda}$ is finite for $\beta\not=\beta'$ by
Corollary~\ref{predcor}~(iii). $\hfill{\square}$

Using the above lemma we prove further properties of $G_H$.
\begin{lemma}\label{desc2}\quad
Let $G_H$ be as above and define $G^\alpha\ (\alpha<\lambda)$ by
$\ G^\alpha:=\{g\in G_H\,|\,||g||<\alpha,\,||g||_{ring}<\alpha\}$.
Then:
\begin{itemize}
\item[{\bf(a)}]
$G_H \cap \widehat{P_\beta}\subseteq G_{\beta+1}$ for all
$\beta<\lambda^*$;
\item[{\bf(b)}]
$\{G^\alpha\,|\, \alpha<\lambda\}$ is a $\lambda$-filtration of
$G_H$; \ and
\item[{\bf(c)}] if $\beta<\lambda^*,\, \alpha<\lambda$ are
ordinals such that $||\varphi_\beta||=\alpha$ then
$G^\alpha\subseteq G_\beta$.
\end{itemize}
\end{lemma}
Note, we used the lower index ($\beta<\lambda^*$) for the
construction while we use the upper index ($\alpha<\lambda$) for
the filtration.\\
{\bf Proof:} First we show (a). Let $g \in G_H \cap
\widehat{P_{\beta}}$ for some $\beta<\lambda^*$. Since
$B_{\beta+1} \subseteq G_{\beta+1}$ we assume $g\in G_H\setminus
B_{\beta+1}$. Then, by Lemma~\ref{desc1}, $g \in B +
\sum\limits_{\gamma\in N}R_{\gamma} y^{(k)}_{\gamma}$ for some
finite $N\subseteq\lambda^*,\,k < \omega$ such that
$[g]_{\lambda}\cap
[_{\lambda}y_\gamma]_{\lambda}$ is infinite for $\gamma\in N$.\\ 
Since $g\in\widehat{P_\beta}$ we also have
$[g]_{\lambda}\subseteq [P_\beta]_{\lambda}\, (=[\widehat{P_\beta}]_{\lambda})$.\\ 
If $||g|| < ||P_\beta||$ then $N\subseteq\beta$ by
Corollary~\ref{predcor}(ii) and thus $g\in
G_\beta\subseteq G_{\beta+1}$.\\ 
Otherwise, if $||g||=||P_\beta||\,(\in \lambda^o)$ then $||g|| =
||y_{\gamma_*}||=||\varphi_{\gamma_*}||$ for $\gamma_*=\max N$ and
$[g]_{\lambda}\cap [y_{\gamma_*}]_{\lambda}\subseteq
[\varphi_\beta]_{\lambda}\cap [\varphi_{\gamma_*}]_{\lambda}$ is
infinite. Hence $\beta=\gamma_*$ by condition (iii) of
Corollary~\ref{predcor} and so $g\in G_{\beta+1}$ as required.
\\ 
Condition (b) is obvious.\\ 
To see (c) let $\beta<\lambda^*,\,\alpha<\lambda$ with
$||\varphi_\beta||=\alpha$ and let $g\in G^\alpha$.  If $g\in
B_{\beta}$ we are finished. Otherwise, by Lemma~\ref{desc1}, we
have $g \in B\, +\, \sum\limits_{\gamma\in N}R_{\gamma}
y^{(k)}_\gamma\ (N\subseteq\lambda^*\text{ finite},\, k\in\omega)$
with $[g]_{\lambda}\cap [y_\gamma]_{\lambda}$ is infinite for
$\gamma\in N$. This implies $||\varphi_\gamma||=||y_\gamma||\leq
||g||<\alpha= ||\varphi_\beta||$ for all $\gamma \in N$ and thus
$N\subseteq\beta$ by Corollary~\ref{predcor}(ii), i.e. $\, g\in
G_\beta$, which finishes the proof. $\hfill{\square}$

\section{Proof of the Main Theorem}

In this final section we want to prove our Main Theorem which
reads as follows:

\begin{mtheorem}
\label{main} Let $\lambda$ be an infinite cardinal such that
$\lambda=\mu^+$ with $\mu^{\aleph_0} = \mu$ and let $H$ be an
$\aleph_1$-free abelian group of size less than $\lambda$. Then
there exists an $\aleph_1$-free $\E(H)$-group of cardinality
$\lambda$ with non-commutative endomorphism ring. In particular,
there is a strong $\E$-group of cardinality $\lambda$.
\end{mtheorem}

{\bf Proof:} Let $\lambda$ and $H$ be given as stated in the
theorem and choose a stationary subset $E$ of $\lambda$ whose
members have cofinality $\omega$. We construct $q_H=(R_H, G_H,
\sigma_H, \sigma_{H*}) =\bigcup\limits_{\alpha <
\lambda}q_{\alpha}$ as in the previous section. Thus $R_H
\subseteq_* \End_{\Z}G_H$ and $\sigma: G_H \rightarrow R_H$ with
kernel $\ker \sigma_H=H$. Moreover, $G_H$ is $\aleph_1$-free and
$R_H$ is obviously non-commutative. We first claim that
$\sigma_H$ is surjective. Therefore let $r \in R_H$, hence there
exists $\alpha < \lambda$ such that $r \in R_{\alpha}$. Without
loss of generality we may assume that $\alpha \not\in E$. By (v)
we conclude that $r \in R_{\alpha} \subseteq \Im \sigma_{\alpha
+1}$. Since $\sigma_{\alpha +1} \subseteq \sigma_H$ as functions
it follows that $r \in \Im \sigma_H$ and thus $\sigma_H$ is
surjective.\\
It remains to prove that $R_H=\End_{\Z}G_H$. Assume that $\varphi
\in \End_{\Z}G_H \backslash R$. Let $\varphi'=\varphi
\restriction_{B'}$, hence $\varphi' \not\in R$. Let $I=\{
(\varepsilon_n, \alpha_n) | n < \omega \} \subseteq \rho \times
\lambda$ such that $\alpha_0 < \alpha_1 < \cdots \alpha_n <
\cdots$ and $I_{\lambda} \cap [g]_{\lambda}$ is finite for all $g
\in G_H$. Note that the existence of $I$ can be easily arranged,
e.g. let $E \subsetneq \lambda^0, \alpha \in \lambda^0 \backslash
E, \varepsilon_n \in \rho (n < \omega)$ arbitrary and
$(\alpha_n)_{n < \omega}$ any ladder on $\alpha$.\\
By Lemma \ref{step} there exists an element $y \in \widehat{B'}$
such that $y \varphi' \not\in G_H'$ which is an extension of
$G_H$. By the Black Box Theorem \ref{black} the set \[ E'=\{
\alpha \in E | \exists \beta < \lambda^* : ||\varphi_{\beta}
||=\alpha, \varphi_{\beta} \subseteq \varphi', [y] \subseteq
[\varphi_{\beta}] \} \] is stationary since $|[y]| < \aleph_0$.
Note, $[y]\subseteq [\varphi_{\beta}]$ implies that $y \in
\widehat{\dom \varphi_{\beta}}$. Moreover, the set $C=\{ \alpha <
\lambda : \varphi(G^{\alpha}) \subseteq G^{\alpha} \}$ is a cub in
$\lambda$, hence $E' \cap C \not= \emptyset$. Let $\alpha \in E'
\cap C$. Then $G^{\alpha} \varphi' \subseteq G^{\alpha}$ and
there exists an ordinal $\beta < \lambda^*$ such that $||
\varphi_{\beta}|| =\alpha$, $\varphi_{\beta} \subseteq \varphi$
and $y \in \widehat{\dom \varphi_{\beta}}$. The first property
implies that $G^{\alpha} \subseteq G_{\beta}$ by Lemma
\ref{desc2} and the latter properties imply that $\varphi_{\beta}
\not\in R$. Moreover, $\dom \varphi_{\beta} \subseteq B'$ with
$|| \dom \varphi_{\beta} ||_{ring} \leq ||\dom
\varphi_{\beta}||=\alpha$ and hence $\dom \varphi_{\beta}$, and
also $(\dom \varphi_{\beta})\varphi$ are contained in $G^{\alpha}
\subseteq G_{\beta}$.\\
Therefore $\varphi_{\beta} : \dom \varphi_{\beta} \rightarrow
G_{\beta}$ with $\varphi_{\beta} \not\in R_{\beta}$ and thus it
follows form the construction that $y_{\beta} \varphi_{\beta}
\not\in G_{\beta +1}$. On the other hand it follows from Lemma
\ref{desc2} that $y_{\beta} \varphi_{\beta} = y_{\beta} \varphi
\in G_H \cap \widehat{P_{\beta}} \subseteq G_{\beta +1}$ - a
contradiction. Thus $\End_{\Z}G_H=R$. $\hfill{\square}$

\begin{corollary}
Let $\lambda$ be an infinite cardinal such that
$\lambda^{\aleph_0}=\lambda$. Then there is a non-commutative
ring $R$ such that $\End_{\Z}(R_{\Z}) \cong R$.
\end{corollary}

{\bf Proof:} Let $H=\{0\}$ and apply Theorem \ref{main} to obtain
a strong $\E$-group $G$. Thus $\sigma: G \rightarrow \End_{\Z}G$
is an isomorphism, where $\End_{\Z}G$ is non-commutative. Hence,
$R:=G$ has a non-commutative ring structure such that
$\End_{\Z}(R_{\Z}) \cong R$.
 $\hfill{\square}$

\noindent
R\"udiger G\"obel and Lutz Str\"ungmann \\
Fachbereich 6, Mathematik und Informatik \\
Universit\"at Essen, 45117 Essen, Germany \\
{\small e--mail: R.Goebel@Uni-Essen.De\\
\hspace*{1.2cm} lutz.struengmann@uni-essen.de}\\
and \\
Saharon Shelah \\
Department of Mathematics\\
Hebrew University, Jerusalem, Israel \\
and Rutgers University, Newbrunswick, NJ, U.S.A \\
{\small e-mail: Shelah@math.huji.ae.il}

\end{document}